\documentclass[12pt,twoside]{article}
\title{On strong n-perfect rings}
\date{}
\usepackage{amsfonts}
\usepackage{amsmath}
\usepackage{amssymb}
\usepackage{latexsym}
\newtheorem{thm}{\bf Theorem}[section]
\newtheorem{cor}[thm]{\bf Corollary}
\newtheorem{lem}[thm]{\bf Lemma}
\newtheorem{prop}[thm]{\bf Proposition}

\newtheorem{rem}[thm]{\bf Remark}

\newtheorem{exmp}[thm]{\bf Example}

\catcode`\ç=13
\defç{\c{c}}
\catcode`\é=13
\defé{\'e}
\catcode`\à=13
\defà{\`a}
\catcode`\è=13
\defè{\`e}
\catcode`\â=13
\defâ{\^a}
\catcode`\ù=13
\defù{\`u}
\catcode`\ê=13
\defê{\^e}
\catcode`\î=13
\defî{\^\i}
\catcode`\ô=13
\defô{\^o}
\newcommand{\field}[1]{\mathbb{#1}}
\newcommand{\C}{\field{C}}

\newcommand{\Q }{\field{Q}}
\newcommand{\Z }{\field{Z}}

\def\proof{{\parindent0pt {\bf Proof.\ }}}

\def\wdim{{\rm wdim}}

\def\gldim{{\rm gldim}}

\def\pd{{\rm pd}}

\def\fd{{\rm fd}}

\def\sup{{\rm sup}}

\def\qf{{\rm qf}}

\def\Ext{{\rm Ext}}

\def\Tor{{\rm Tor}}

\newcommand{\cqfd}
{\hspace{1cm}
\rule{2mm}{2mm}%
\medbreak%
\par%
}
\begin{document}
\thispagestyle{empty}
%%%%%%%%%%%%%%%%%%%%%%%%%%%%%%%%%%%%%%%%%%%%%%%%%%%%%%%%%
%%%%%%%%%%%%%%%%%%%%%%%%%%%%%%%%%%%%%%%%%%%%%%%%%%%%%%%%%
%%%%%%%%%%%%%%%%%%%%%%%%%%%%%%%%%%%%%%%%%%%%%%%%%%%%%%%%%
%%%TITLE%%%%%%%%%%%%%%%%%%%%%%%%%%%%%%%%%%%%%%%%%%%%%%%%%
\maketitle \vspace*{-2cm}
\begin{center}{\large\bf Abdellatif Jhilal and Najib Mahdou}
%%%%%%%%%%%%%%%%%%%%%%%%%%%%%%%%%%%%%%%%%%%%%%%%%%%%%%%%%
%%%%%%%%%%%%%%%%%%%%%%%%%%%%%%%%%%%%%%%%%%%%%%%%%%%%%%%%%
%%%%%%%%%%%%%%%%%%%%%%%%%%%%%%%%%%%%%%%%%%%%%%%%%%%%%%%%%
%%%NAMES%%%%%%%%%%%%%%%%%%%%%%%%%%%%%%%%%%%%%%%%%%%%%%%%%

\bigskip
%%%%%%%%%%%%%%%%%%%%%%%%%%%%%%%%%%%%%%%%%%%%%%%%%%%%%%%%%
%%%%%%%%%%%%%%%%%%%%%%%%%%%%%%%%%%%%%%%%%%%%%%%%%%%%%%%%%
%%%%%%%%%%%%%%%%%%%%%%%%%%%%%%%%%%%%%%%%%%%%%%%%%%%%%%%%%
%%%%%%%%%%%%ADDRESSES%%%%%%%%%%%%%%%%%%%%%%%%%%%%%%%%%%%%%%%%%%%%%
\small{Department of Mathematics, Faculty of Science and Technology
of Fez,\\ Box 2202, University S. M.
Ben Abdellah Fez, Morocco, \\  Jhilalabdo@hotmail.com \\
mahdou@hotmail.com}
\end{center}

\bigskip\bigskip
%%%%%%%%%%%%%%%%%%%%%%%%%%%%%%%%%%%%%%%%%%%%%%%%%%%%%%%%%
%%%%%%%%%%%%%%%%%%%%%%%%%%%%%%%%%%%%%%%%%%%%%%%%%%%%%%%%%
%%%%%%%%%%%%%%%%%%%%%%%%%%%%%%%%%%%%%%%%%%%%%%%%%%%%%%%%%
%%%ABSTRACT%%%%%%%%%%%%%%%%%%%%%%%%%%%%%%%%%%%%%%%%%%%%%%
\noindent{\large\bf Abstract.}\bigskip

 In this paper we introduce the notion of ``strong $n$-perfect rings" which is in some way a
generalization of the notion of ``$n$-perfect rings". We are mainly
concerned with those class of rings in the context of pullbacks.
Also we exhibit a class of $n$-perfect rings that are not strong
$n$-perfect rings. Finally, we establish the transfer of this notion
to the direct product notions.\\
\bigskip
%%%%%%%%%%%%%%%%%%%%%%%%%%%%%%%%%%%%%%%%%%%%%%%%%%%%%%%%%

\small{\noindent{\bf Key Words.} Perfect ring, $n$-perfect ring,
strong $n$-perfect property, $(n,d)$-property,
 pullback ring.

%%%%%%%%%%%%%%%%%%%%%%%%%%%%%%%%%%%%%%%%%%%%%%%%%%%%%%%%%
%%%%%%%%%%%%%%%%%%%%%%%%%%%%%%%%%%%%%%%%%%%%%%%%%%%%%%%%%
%%%INTRODUCTION%%%%%%%%%%%%%%%%%%%%%%%%%%%%%%%%%%%%%%%%%%
\begin{section}{Introduction}\bigskip
Throughout this work, all rings are commutative with identity
element, and all modules are unitary. Let $R$ be a ring and let $M$
be an $R$-module. As usual we use $\pd_R(M)$ and $\fd_R(M)$ to
denote the usual projective and flat dimensions of $M$,
respectively. $\gldim(R)$ and $\wdim(R)$ are, respectively , the
classical global and weak dimension of $R$. If $R$ is an integral
domain, we denote its quotient field by $\qf(R)$.\bigskip

A ring $R$ is perfect if every flat $R$-module is  projective
$R$-module. The pioneering work on perfect rings was done by Bass
\cite{B} and most of the  principal characterizations of perfect
rings are contained in Theorem P from that paper.\bigskip

Recently, Enochs, et al. in his work \cite{En} extended the notion
of perfect rings to $n$-perfect rings, such that a ring is called
$n$-perfect, if every flat module has projective dimension less or
equal than $n$ \cite[Definition 1.1]{En}.\bigskip

The classes of rings we will define here are in some ways
generalization of the notion of ``$n$-perfect rings".  Let $n$ be a
positive integer. $A$ commutative ring $R$ is called a strong
$n$-perfect ring if any $R$-module of flat dimension less or equal
than $n$ have projective dimension less or equal than $n$. Note that
if $n=0$ then the strong $0$-perfect rings are the perfect rings.
And it is trivial to remark that every strong $n$-perfect ring is
also an $n$-perfect ring.\bigskip

 In this paper, we investigate the transfer of strong $n$-perfect property
 in some know ring constructions. Namely, we study the strong
 $n$-perfect property of pullbacks (Theorem \ref{5}) and of finite
 direct product (Theorem \ref{26}). See \cite{Fo} for a quite general development of properties of pullbacks.
  And Theorem \ref{5} enable us to give a class of strong $n$-perfect
 rings (Example \ref{24}). \bigskip

 We also show that, in a particular cases, the relationship between
  the strong $n$-perfect property and $(n,d)$-property  ( Theorem
 \ref{15}), for background on $(n,d)$-rings we refer the reader to \cite{Co, Ck, DKM, M}. \bigskip

 As the strong $n$-perfect property is a general of the $n$-perfect
 property, it is important to give examples of $n$-perfect rings
 which are not strong $n$-perfect rings. This is  given in Example
 \ref{25}. \bigskip

 Finally, the Proposition \ref{13} allow us to study the strong $n$-perfect property of some
  particular rings.

\bigskip

\end{section}
%%%%%%%%%%%%%%%%%%%%%%%%%%%%%%%%%%%%%%%%%%%%%%%%%%%%%%%%
%%%%%%%%%%%%%%%%%%%%%%%%%%%%%%%%%%%%%%%%%%%%%%%%%%%%%%%%%%%
%%%%%%%%%%%%%%%%%%%%%%%%%%%%%%%%%%%%%%%%%%%%%%%%%%%%%%%%
%%%%%%%%%%%%%%%%%%%%Section 2.%%%%%%%%%%%%%%%%%%%%%%%%%%%%
\begin{section}{Main results}\bigskip

 We begin  by the following change of rings result
  for strong $n$-perfect property, in which we use the notion of
  flat epimorphism of rings, which is defined as follows: Let
 $\Phi:A\rightarrow B$ be a ring homomorphism. $B$ $($ or $\Phi)$ is
 called  a flat epimorphism of $A$,  if $B$ is a flat $A$-module and
 $\Phi$ is an epimorphism, that is, for any two ring homomorphism
 $B\overset{f}{\underset{g}{\rightrightarrows}}C$, $C$ a ring,
 satisfying $f\circ \Phi= g\circ\Phi$, we have $f=g$ \cite[pages
 13-14]{G}.
  For example $S^{-1}A$ is a flat epimorphism of $A$ for every
  multiplicative set $S$ of $A$. Also, the quotient ring $A/I$ is a
  flat epimorphism of $A$ for every pure ideal $I$ of $A$, that is,
  $A/I$ is a flat $A$-module \cite[Theorem 1. 2. 15]{G}

\begin{prop}
Let $A$ and $B$ be two rings such that  $\Phi: A\rightarrow B$ be  a
flat epimorphism of  $A$. Then:
\\ If $A$ is a strong n-perfect ring then $B$  is a strong n-perfect ring.
\end{prop}
\proof: Let $M$ be an $B$-module such that $\fd_{B}(M)\leq n$. Our
aim is to  show that $\pd_{B}(M)\leq n$.
\\ By hypothesis we have $\Tor^{A}_{k}(M,B)=0$ for all $k>0$, so for
any $B$-module $N$, we have from \cite[Proposition 4.1.3]{CE}.
$$(\ast) \qquad  \Ext_{A}^{n+1}(M,N\otimes_{A}B)\cong
\Ext_{B}^{n+1}(M\otimes_{A}B,N\otimes_{A}B)$$ And \cite[Theorem
1.2.19]{G} give $M\otimes_{A}B\cong M$ and $N\otimes_{A}B\cong N$.
\\On the other hand $\fd_{A}(M)\leq \fd_{B}(M)\leq n $
\cite[Exercise 10.p.123]{CE}. Thus $\pd_{A}(M)\leq n$ since $A$ is a
strong n-perfect ring.  Therefore $\Ext_{B}^{n+1}(M,N)=0$ by
$(\ast)$ and so $\pd_{B}(M)\leq n$. \cqfd
\bigskip

 As $S^{-1}A$ and  $A/I$  are the particular cases of flat
 epimorphism of $A$  for any multiplicative set $S$ of $A$ and
 for any pure ideal $I$ of $A$. we have the
  following results.

\begin{cor}\label{8}
\begin{enumerate}
\item For every multiplicative set $S$ of a strong n-perfect ring $A$,
$S^{-1}A$ is a strong n-perfect ring $(n \geq 0)$.
\item If $A$ is a strong n-perfect ring, then so is the quotient
ring $A/I$ for every pure ideal $I$ of $A$ $(n \geq 0)$.

\end{enumerate}
\end{cor}\bigskip

 $\quad$ Now we give our main result of this paper.
 
 \begin{thm}\label{5}
Let $A\hookrightarrow B$ be an injective flat ring homomorphism and
let $Q$ be a pure ideal of $A$ such that $QB=Q$. Then:
\begin{enumerate}
\item If $B$ is a strong n-perfect ring. Then: \\ $A/Q$ is a strong
n-perfect ring if and only if $A$ is a strong n-perfect ring.
\item If $B=S^{-1}A$, where $S$ is a multiplicative set of $A$. Then:
\\ $A$ is a strong $n$-perfect ring if and only if $B$  and $ A/Q$ are a strong
$n$-perfect rings.
\end{enumerate}
\end{thm}\bigskip

Before proving this theorem, we establish the following lemmas.

 \begin{lem}\cite[Lemma 2.5 ]{DKM}\label{4}
Let $A\hookrightarrow B$ be an injective flat ring homomorphism and
let $Q$ be an ideal of $A$ such that $QB=Q$. Let $E$ be an
$A$-module such that $E\otimes_{A}B$ is $B$-flat.Then:\\
$\pd_{A}(E)\leq d \Leftrightarrow \pd_{B}(E\otimes_{A}B)\leq d  $
and $ \pd _{A/Q}(E\otimes_{A}A/Q)\leq d$

 \end{lem}\bigskip

\begin{lem}\label{6}Let $A\hookrightarrow B$ be an injective flat ring homomorphism and
let $Q$ be a pure ideal of $A$ such that $QB=Q$. Let $E$ be an
$A$-module. Then:

\begin{enumerate}
\item $\fd_{A}(E)\leq d \Leftrightarrow
\fd_{B}(E\otimes_{A}B)\leq d  $ and $ \fd
_{A/Q}(E\otimes_{A}A/Q)\leq d$
\item $\pd_{A}(E)\leq d \Leftrightarrow \pd_{B}(E\otimes_{A}B)\leq d $
and $ \pd _{A/Q}(E\otimes_{A}A/Q)\leq d$
\end{enumerate}
\end{lem}\bigskip
\proof Similar to the proof  of Lemma \ref{4}.\cqfd

\textbf{Proof of theorem} \ref{5}:
\\$1)$ If $A$ is a strong $n$-perfect ring, then so is $A/Q$ by
Corollary \ref{8}$(2)$ since $Q$ is a pure ideal of $A$. Conversely,
assume that $B$ and $A/Q$ are a strong $n$-perfect rings. Let $M$ be
an $A$-module such that $\fd_{A}(M)\leq n$. Then
$\fd_{B}(M\otimes_{A}B)\leq n $ and $ \fd
_{A/Q}(M\otimes_{A}A/Q)\leq n$ by lemma \ref{6}.1. Hence,
$\pd_{B}(M\otimes_{A}B)\leq n  $ and $ \pd
_{A/Q}(M\otimes_{A}A/Q)\leq n$, since $B$ and $A/Q$ are a strong
n-perfect rings. So that $\pd_{A}(M)\leq n$  by lemma \ref{6}.(2),
therefore $A$ is a strong n-perfect ring.
\\$2)$ By Corollary \ref{8} and $1)$ and this completes the proof of
Theorem \ref{5}. \cqfd  \bigskip

From this  Theorem we deduce easily the following  example of strong
$n$-perfect rings.
\begin{exmp}\label{24}
 Let $D$ be an integral domain such that $\gldim(D)=n$, let
$K=\qf(D)$ and let $n\geq 2$.
  \\ Consider the quotient ring  $S=K[X]/(X^{n}-X)=K+ \overline{X}K[\overline{X}] = K +
  I$ with $I= \overline{X}K[\overline{X}]$.
  Set $R=D + I$.
   Then $R$ is a strong $n$-perfect ring.

\end{exmp}
\proof First we show that $I$ is a pure ideal of $R$. Let $\alpha$
be an element of $I$, then $\alpha= \overline{X}^{i} ( a_{0} +
a_{1}\overline{X}+ ... + a_{n-1}\overline{X}^{n-1})$ with $a_{i}\in
K$ for $1\leq i \leq n-1$, and $a_{0}\neq 0$. As $\delta= a_{0} +
a_{1}\overline{X}+ ... + a_{n-1}\overline{X}^{n-1}$ is an invertible
element.\\
 And  $\overline{X}^{i}(1-\overline{X}^{n-1})=\overline{X}^{i}-\overline{X}^{n+(i-1)}=\overline{X}^{i}-
\overline{X}^{n}\overline{X}^{i-1}=\overline{X}^{i}-\overline{X}^{i}=\overline{0}$

Then $\alpha ( 1- \overline{X}^{i-1})= 0$ and $\overline{X}^{i-1}\in
I$. Therefore $I$ is a pure ideal of $R$ by \cite[Theorem
1.2.15]{G}.
\\ On the other hand, we have $S$ an artinian ring. According to \cite[Corollary 28.8]{A}
we deduce $S$ is a perfect ring and so $S$ is a strong $n$-perfect
ring, and also $D$ is a strong $n$-perfect ring, since
$\gldim(D)=n$. Therefore $R$ is a strong $n$-perfect ring by Theorem
\ref{5} . \cqfd
\bigskip

Later, we give an example showing that Theorem \ref{5} is not true
without assuming that the ideal $Q$ is pure (see Example
\ref{7}).\bigskip

Now we link the strong $n$-perfect property with $(n,d)$-property in
a particular cases. First, let $n,$ $d$ be positive  integers, we
say that $R$ is an $(n, d)$-ring, if each $n$-presented $R$-module
has projective
dimension at most $d$, see for instance \cite{Co, Ck, DKM, M}.\\
In particular, the $(1,0)$-rings are the von Neumann regular rings,
$(0,1)$-rings are the hereditary rings, $(1,1)$-rings are the
semihereditary rings and $(0,0)-$rings are the semisimple rings
\cite[Theorem 1.3]{Co}.\bigskip

It is clear to see that if $R$ is Noetherian, then $R$ is an $(1,
d)$-ring if and only if $R$ is an $(0, d)$-ring. In the context of
strongly $d$-perfect ring, we have:

\begin{thm}\label{15}
Let $R$ be a ring.
\\$R$ is an $(1,d)$-ring and a strong $d$-perfect ring if and only if $R$
is an $(0,d)$-ring.
\end{thm}
\proof:  Assume that $R$ is an $(1,d)$-ring and strong $d$-perfect.
Hence $\wdim(R)\leq d$ \cite [Theorem 1.3.9]{G} since for finitely
generated $I$ ideal of $R$, we have $\fd(R/I)\leq d$ (since
$\wdim(R)\leq d$).
\\ New let $M$ be an $R$-module, thus $\fd(M)\leq d$ (since $\wdim(R)\leq d$) and so $\pd(M)\leq d$ since
 $R$ is a
strong $d$-perfect ring . Therefore $R$ is an $(0,d)$-ring.
\\ Conversely, assume that  $R$ is an $(0,d)$-ring that is $\gldim(R)\leq d$. It is clear
 that $R$ is an
$(1,d)$-ring and  a strong d-perfect ring. \cqfd
 \bigskip

According to proof of this theorem, we also show that strong
$n$-perfect property links, in particular cases, the classical
global and weak dimensions.

\begin{cor}\label{11}
Let $R$ be a strong $n$-perfect ring then: $\gldim(R)\leq n$ if and
only if $\wdim(R)\leq n$.

\end{cor}
\bigskip

 Since
the  $(1,0)$-rings are the von Neumann regular rings, the
$(0,0)-$rings are the semisimple rings, the $(0,1)$-rings are the
hereditary rings, and the  $(1,1)$-rings are the semihereditary
rings, then by Theorem \ref{15} we have:
\begin{cor}\label{9}
\begin{enumerate}
\item A von Neumann regular ring is a semisimple ring if and only if
is a perfect ring.
\item A semihereditary ring is an hereditary ring if and only if is a strong
$1$-perfect ring.

\end{enumerate}
\end{cor}

\begin{rem}
There exists a ring $R$ which is a von Neumann regular ring and an
hereditary ring but is not a semisimple ring \cite[Example 2.7]{Co}.
By corollary \ref{9}, we deduce that $R$ is a strong $1$-perfect
ring but it is not a perfect ring.

\end{rem}

 Now we can  show that Theorem
\ref{5} is not true without assuming that $Q$ is a pure ideal.
\begin{exmp}\label{7}  Let $\Z$ denotes the ring of integers, and let $\Q$
denotes the field of rational numbers. \\ Consider the discrete
valuation domain $B=\Q[X]_{(X)}=\Q + XB$. Let $A=\Z + XB$, then.
\begin{enumerate}
\item $B$ is a strong 1-perfect ring.
\item $A/XB\cong \Z$ is a strong 1-perfect ring.
\item $XB$ is not a pure ideal, and $A$ is not a strong 1-perfect ring.

\end{enumerate}

\end{exmp}

\proof: 1) We have $B$ is a strong 1-perfect ring since $B$ is an
hereditary ring (corollary \ref{9}).
\\2) Similarly $A/XB\cong\Z$ is a strong 1-perfect ring.
\\3) We first set $M=XB$ and $E=A/M$ thus $A=\Z+M$ and $B=\Q+M$.
The exact sequence of $A$-modules  $0\rightarrow M\rightarrow
A\rightarrow \Z\rightarrow 0$ yields the exact sequence of
$A$-modules
$$0\rightarrow \Tor_{1}^{A}(E,\Z)\rightarrow E\otimes_{A}M\rightarrow
E\rightarrow E/ME \rightarrow 0$$   But $\Tor_{1}^{A}(E,\Z)\cong
E\otimes_{A}M\cong M/M^{2} \neq 0$ Since $M$ is principal and
$ME=0$, i.e., $\Tor^{A}_{1}(A/M, \Z)\neq 0$. So $M$ is not a pure
ideal \cite [Theorem 1.2.15]{G}. \\ On the other hand: we have
$\fd_{A}(E)\leq 1$, since $0\rightarrow M\rightarrow A\rightarrow
\Z(E)\rightarrow 0$ is exact and $M$ is $B$-plat and $B$ is
$A$-plat.
\\If $\pd_{A}(E)\leq 1$ was true, we would have an $A$-projective
resolution: $$0\rightarrow P_{1}\rightarrow
P_{0}\stackrel{u_{0}}\rightarrow E\rightarrow 0$$

Thus $$(\ast) \qquad 0\rightarrow \Tor_{1}^{A}(E,\Z)\rightarrow
P_{1}\otimes_{A}\Z\rightarrow P_{0}\otimes_{A}\Z
\stackrel{u_{0}\otimes1}\rightarrow E\otimes_{A}\Z\rightarrow 0$$ is
exact. But $E\otimes_{A}\Z\cong E\otimes _{A}A/M\cong E/ME\cong E=
\Z$, thus $E\otimes_{A}\Z$ is $\Z$-projective, making $(\ast)$
split-exact and giving that $P_{1}\otimes_{A}\Z\cong
\Tor_{1}^{A}(E,\Z)\oplus Ker(u_{0}\otimes1)$, so that
$\Tor_{1}^{A}(E,\Z)\cong M/M^{2}$ is a projective $\Z$-module. But
$M/M^{2}$ is a $\Q$-vector space, so it would follow that $\Q$ is a
projective $\Z$-module, which is false since $\Q=qf(\Z)$ and $\Z\neq
\Q$. Thus $\pd_{A}(E)>1$. \cqfd \bigskip

But the transfer of the ``$n$-perfectness property'' in pullback
rings
 is not need to the hypothesis ``$Q$ is a pure ideal'' as follows to explain.

\begin{prop}\label{10}
Let $A\hookrightarrow B$ be an injective flat ring homomorphism and
let $Q$ be an ideal of $A$ such that $QB=Q$. Then:
\\If $B$ is an n-perfect ring and $A/Q$ is an n-perfect ring, then $A$ is an n-perfect ring.
\end{prop}
\proof Let $E$ be a flat $A$-module  then $E\otimes_{A}B$ is a flat
$B$-module. \\Thus $\pd_{B}(E\otimes_{A}B)\leq n$ since $B$ is an
n-perfect ring.\\ On the other hand , $E\otimes_{A}A/Q$ is a flat
$A/Q$-module, thus $\pd _{A/Q}(E\otimes_{A}A/Q)\leq n$ since $A/Q$
is
 an n-perfect ring. By lemma \ref{4}, we have $\pd_{R}(E)\leq n$. Therefor $A$ is
an n-perfect ring. \cqfd \bigskip

In example \ref{7}, the ring $A$ is an $1$-perfect ring such that is
not a strong $1$-perfect ring. The following example shows that
there exists the rings $n$-perfect are not a strong $n$-perfect
rings for any $n$.\\
\begin{exmp}\label{25}
Let $\C$ denotes the field of complex numbers and $\C(X,Y)$ denotes
the quotient field of the polynomial ring $\C[X,Y]$. Consider the
discrete valuation domain $V=\C(X,Y)[Z]_{(Z)}=\C(X,Y)+ZV$.
\\ Set the pullback ring $R=\C[X,Y]+ZV$. Then:
\begin{enumerate}
\item $R$ is an $2$-perfect ring and it is not a strong $2$-perfect ring .
\item $R[X_{1}]$ is an $3$-perfect ring and it is not a strong  $3$-perfect.
\item Let $n\geq 4$ be integer $R[X_{1},...,X_{n-2}]$ is an $n$-perfect ring and it is
not strong $n$-perfect ring.

\end{enumerate}
\end{exmp}
 \proof: By \cite [Proposition 2.1]{Do} $\gldim(R)=3$ and
 $\wdim(R)=2$.
 \\ $1)$ We have
 $R/ZV= \C[X,Y]$ is an $2$-perfect ring (since $\gldim(\C[X,Y])=\wdim(\C[X,Y])=2$)
  and   $V$ is an $2$-perfect ring (since $V$ is a
discrete valuation domain).  Hence, $R$ is an $2$-perfect ring (By
proposition \ref{10}) and  $R$ is not a strong $2$-perfect  (by
Corollary \ref{11}).
\\ $2)$ $R[X_{1}]$ is an $3$-perfect ring (by  \cite [Example iv page
2]{En}) and  $R$ is not
 strong $3$-perfect (by  $\wdim(R[X_{1}])=3$ and $\gldim(R[X_{1}])=4$).
 \\3) By induction on $n$, we have for every $n\geq 4$,
 $R[X_{1},...,X_{n-2}]$ is an  $n$-perfect ring and it is not  a strong
 $n$-perfect ring.\cqfd \bigskip

In the following we prove that the strong $n$-perfect property
descends into a faithfully flat ring homomorphism.

\begin{prop}\label{13} Let $A\longrightarrow B$ be a ring homomorphism making
$B$ a faithfully flat $A$-module. if $B$ is a strong $n$-perfect
ring then $A$ is a strong $n$-perfect ring.

\end{prop}
\proof: Let $M$ be an $A$-module such that $\fd_{A}(M)\leq n$. Our
aim is to show that $\pd_{A}\leq n$.
\\ Let $0\longrightarrow P \longrightarrow F_{n-1}\longrightarrow
...\longrightarrow F_{1}\longrightarrow F_{0}\longrightarrow M
\longrightarrow 0$ be an exact sequence of $A$-modules, where
$F_{i}$ is a free $A$-module for each $i$, $P$ is a flat $A$-module
thus:
\\ $(\ast)$ $0\longrightarrow P\otimes_{A}B \longrightarrow F_{n-1}\otimes_{A}B
\longrightarrow ...\longrightarrow F_{1}\otimes_{A}B\longrightarrow
F_{0}\otimes_{A}B\longrightarrow M\otimes_{A}B \longrightarrow 0$
 is an exact sequence of $B$-modules, where $F_{i}\otimes_{A}B$ is a
 free $B$-module for each $i$ and $P\otimes_{A}B$ is a flat $B$-module,
 therefore $\fd_{B}(M\otimes_{A}B)\leq n$ and so  $\pd_{B}(M\otimes_{A}B)\leq n$ $B$ is a strong
 $n$-perfect. Hence,
 $P\otimes_{A}B$ is a projective $B$-module   By $(\ast)$. \\According to
 \cite[Example 3.1.4. page 82]{Gr} we deduce that $P$ is a projective $A$-module
 since $P$ is a flat $A$-module. Then $\pd_{A}(M)\leq n$ and so
  $A$ is a strong $n$-perfect ring.\cqfd \bigskip
  We profit by this proposition to study the strong $n$-perfect
property of some particular rings.\bigskip

 Recall that the trivial
extension of a ring $R$ by an $R$-module $M$  is the ring denotes by
$R\alpha M$ whose underling group is $R\times M$ with pairwise
addition and multiplication given by $(a, e)(a',e')=(aa',ae'+a'e)$
(see for instance \cite{F,G,H}).
\begin{cor}\label{19}
\begin{enumerate}
\item Let $A$ be a ring, let $M$ be a flat $A$-module and let $R= A\alpha
M$ be the trivial ring extension of $A$ by $M$. If $R$ is a strong
$n$-perfect ring then  so is $A$ .
\item Let $A\subset B$ be two ring such that $B$ is a flat
$A$-module. Let $S=A+XB[X]$, where X is an indeterminate over $B$ .
If $S$ is a strong $n$-perfect ring then  so is $A$ .
\item Let $R$ be a ring and $X$ is an indeterminate over $R$. If $R[X]$ is
 a strong \\$n$-perfect ring  then so is $R$ .
\end{enumerate}

\end{cor}
\proof: 1) We have $M$ a flat $A$-module, thus the ring homomorphism
$ A\hookrightarrow R$ is a faithfully flat, by Proposition \ref{13},
we conclude that $A$ is a strong $n$-perfect ring since $R$ is a
strong $n$-perfect ring.
\\ 2) The ring $B$ is a flat $A$-module and $XB[X]\cong B[X]$ thus
$S= A + XB[X]$ is a faithfully flat $A$-module, also by Proposition
\ref{13} the ring $A$ is a  strong $n$-perfect  since $S$ is a
strong $n$-perfect ring.
\\ 3) Similarly, by Proposition \ref{13}, since $R[X]$  is a faithfully flat
$R$-module.\cqfd \bigskip

Next  we establish the transfer of
the strong $n$-perfect property to finite direct product.\\

\begin{thm}\label{26} Let $(R_{i})_{i=1,...,m}$ be a family of rings. Then
$\prod_{i=1}^{m}R_{i}$ is a strong n-perfect ring if and only if
$R_{i}$ is a strong n-perfect ring for each $i=1,...,m$ .

\end{thm}
  \quad The proof of the theorem involves the following results.
 \begin{lem}\cite[Lemma 2.5.(2)]{M}\label{2}  \\Let $(R_{i})_{i=1,2}$ be a family of rings and $E_{i}$
 be an $R_{i}$-module for $i=1,2$. We have : $\pd_{R_{1}\times R_{2}}(E_{1}\times
 E_{2})=$\sup$\{\pd_{R_{1}}(E_{1}),\pd_{R_{2}}(E_{2})\}$.

 \end{lem}

 \begin{lem} \label{3} Let $(R_{i})_{i=1,2}$ be a family of rings and $E_{i}$
 be an $R_{i}$-module for $i=1,2$. Then $\fd_{R_{1}\times R_{2}}(E_{1}\times
 E_{2})=$\sup$\{\fd_{R_{1}}(E_{1}),\fd_{R_{2}}(E_{2})\}$.
 \end{lem}
 \proof The proof  is analogous to the proof of Lemma
 \ref{2}.\\ \\
 \textbf{Proof of Theorem \ref{26}}\\ By induction on m, it suffices to prove the
 assertion for $m=2$. Let $R_{1}$ and $R_{2}$ be two rings such that
 $R_{1}\times R_{2}$ is a strong n-perfect
 and let $E_{1}$ be an $R_{1}$-module such that $\fd_{R_{1}}(E_{1})\leq n$,
 and let $E_{2}$ be an $R_{2}$-module such that $\fd_{R_{2}}(E_{2})\leq
 n$. Then
$\fd_{R_{1}\times R_{2}}(E_{1}\times
 E_{2})=$\sup$\{\fd_{R_{1}}(E_{1}),\fd_{R_{2}}(E_{2})\}$ by lemma \ref{3}.
 \\Thus $\pd_{R_{1}\times R_{2}}(E_{1}\times
 E_{2})\leq n \quad$   since   $ \quad R_{1}\times R_{2}$ is a strong n-perfect ring.
 \\But $\pd_{R_{1}\times R_{2}}(E_{1}\times
 E_{2})=$\sup$\{\pd_{R_{1}}(E_{1}),\pd_{R_{2}}(E_{2})\}$ by lemma \ref{2}. \\Therefore
  $\pd_{R_{1}}(E_{1})\leq n$ and
$\pd_{R_{2}}(E_{2})\leq n$ and so $R_{1}$ and $R_{2}$ are  a strong
n-perfect rings .
\\Conversely, assume that $R_{1}$ and $ R_{2}$ be two strong n-perfect
rings. Let $E_{1}\times E_{2}$ be an $R_{1}\times R_{2}$-module
where $E_{i}$ is an $R_{i} $-module for each $i=1,2 ,$ such that
$\fd_{R_{1}\times R_{2}}(E_{1}\times
 E_{2})\leq n. $
 Thus
 $\fd_{R_{1}}(E_{1})\leq n$ and
$\fd_{R_{2}}(E_{2})\leq n$ (by lemma \ref{3})  and so
$\pd_{R_{1}}(E_{1})\leq n$ and $\pd_{R_{2}}(E_{2})\leq n$ (since
$R_{1}$ and $R_{2}$ are strong n-perfect rings). Therefore
$\pd_{R_{1}\times R_{2}}(E_{1}\times
 E_{2})\leq n $ (By Lemma \ref{2}) and so $R_{1}\times R_{2}$ is a strong n-perfect
 rings.\cqfd \bigskip

Recall that a ring $R$ is called an $S$-ring if every finitely
generated flat $R$-module is projective (see \cite{P}).\bigskip

The object of this remark is to conclude, in general case, that the
strong $n$-perfect property is not characterized by the cyclic
modules.

\begin{rem}
The ring $A$ of Example \ref{7} is an integral domain, hence, by
\cite[page 288]{E}, the ring $A$ is an $S$-ring, and by
\cite[Corollary 1.7]{V} we deduce on the ring $A$ every cyclic flat
module is projective. But according also to Example \ref{7} the ring
$A$ is not perfect since it is not strong $1$-perfect. Therefore
 the perfect property is not characterized by the cyclic modules.

\end{rem}

\end{section}

%%%%%%%%%%%%%%%%%%%%%%%%%%%%%%%%%%%%%%%%%%%%%%%%%%%%%%%%%
%%%%%%%%%%%%%%%%%%%%%%%%%%%%%%%%%%%%%%%%%%%%%%%%%%%%%%%%%
%%%REFERENCES%%%%%%%%%%%%%%%%%%%%%%%%%%%%%%%%%%%%%%%%%%%%
%%%%%%%%%%%%%%%%%%%%%%%%%%%%%%%%%%%%%%%%%%%%%%%%%%%%%%%%

%%%%%%%%%%%%%%%%%%%%%%%%%%%%%%%%%%%%%%%%%%%%%%%%%%%%

%%%%%%%%%%%%%%%%%%%%%%%%%%%%%%%%%%%%%%%%%%%%%%%%%%%%

\bigskip\bigskip

%%%%%%%%%%%%%%%%%%%%%%%%%%%%%%%%%%%%%%%%%%%%%%%%%%%%%%%%
\end{document}